\documentclass[11pt,twoside]{amsart}
\usepackage[all]{xy}
        \usepackage {amssymb,latexsym,amsthm,amsmath}
        \usepackage{color}
        \usepackage[lined,boxed,commentsnumbered]{algorithm2e}
\usepackage[pagebackref=true, colorlinks, linkcolor=blue, citecolor=red]{hyperref}
\usepackage{multicol}
        \usepackage{amsfonts}
        \usepackage{a4}
        \usepackage{tikz}
        \usepackage{verbatim}
        \usepackage{caption}
        \usepackage{float}
        
\long\def\symbolfootnote[#1]#2{\begingroup%
\def\thefootnote{\fnsymbol{footnote}}\footnote[#1]{#2}\endgroup}

\newcommand{\n}{\mathfrak n}

\newcommand{\GL}{\textup{GL}}

\makeatletter
\def\imod#1{\allowbreak\mkern10mu({\operator@font mod}\,\,#1)}
\makeatother

\newtheorem{theorem}{Theorem}[section]
\newtheorem{theoremalph}{Theorem}[section]

\newtheorem{corollary}[theorem]{Corollary}

\newtheorem*{theorem*}{Theorem}
\theoremstyle{definition}

\newtheorem{example}[theorem]{Example}
\numberwithin{equation}{section}

\newcommand{\ignore}[1]{}

\newcommand{\mynote}[1]{}
\begin{document}
\setcounter{section}{0}
\title{Representation dimension of some finite groups}
\author{Gurleen Kaur}
\address{Indian Institute of Technology Ropar, Punjab 140001 India}
\email{gurleenkaur992gk@gmail.com}
\author{Amit Kulshrestha}
\address{Indian Institute of Science Education and Research Mohali, Knowledge City, Sector 81, Mohali 140 306, India}
\email{amitk@iisermohali.ac.in}
\author{Anupam Singh}
\address{Indian Institute of Science Education and Research Pune, Dr. Homi Bhabha Road, Pashan, Pune 411 008, India}
\email{anupamk18@gmail.com}

\date{}

\begin{abstract}
For a finite group $G$, the representation dimension is the minimal degree of 
a complex faithful representation of $G$. In this article, we compute representation dimension for some finite $p$-groups, in particular, for all $p$-groups of order up to $p^{5}$, their direct products, and groups with certain conditions on nonlinear irreducible characters. We also compute the minimal degree of a 
complex faithful irreducible representation of such groups, if one exists. In the appendix, we present GAP codes to compute these numbers.
\end{abstract}
\maketitle

\section{Introduction}
\noindent Throughout this paper, we restrict our attention to finite groups and complex representations.
Let $G$ be a finite group.
One of the classical and challenging problems in the representation theory of finite groups is to determine the smallest positive integer $n$ such that $G$ embeds into $\operatorname{GL}(n,\mathbb{C})$, i.e., $G$ is isomorphic to a linear group of degree $n$ (over the complex field $\mathbb{C}$). The study of the linear groups has been of natural interest with intensive development by Flannery, Green, James, Jansen, Lusztig, O'Brien, Steinberg, Tiep and Zalesskii (see \cite{BF, Green, James, Jansen, L, Steinberg, TZ}). Since a faithful representation is an isomorphism of $G$ to a subgroup of $\operatorname{GL}(n,\mathbb{C})$, it follows that a finite group is isomorphic to a linear group of degree $n$, if and only if, it has a faithful representation of degree $n$.
The \textit{representation dimension} $\delta(G)$ of a finite group $G$, coined by Cernele et al. \cite{MR2818953}, is defined as follows:
$$\delta(G) := \min\{\operatorname{deg}(\rho) : \rho \text{ is a complex faithful representation of } G\},$$
where $\operatorname{deg}(\rho)$ denotes the degree of $\rho$.                                  Representation dimension has been extensively studied for finite groups of Lie type (see for example L\"{u}beck \cite{Lubeck}, Tiep and Zalesskii~\cite{TZ}), and for finite $p$-groups (see for example Bardestani et al.~\cite{BMS}, 
 Di Martino and Tamburuni~\cite{MR1153979}, Karpenko and Merkurjev ~\cite{KM}).
Another related number $\delta_{irr}(G)$, as defined below, is studied:
$$\delta_{irr}(G) := {\rm min}\{ \operatorname{deg}(\rho) : \rho \text{ is a complex faithful irreducible representation of } G\}.$$

\noindent We note that $\delta_{irr}(G)$ does not exist for all groups. Due to Theorem 2.32 of \cite{Isaacs}, if $\delta_{irr}(G)$ exists, then $\mathcal{Z}(G)$ is cyclic, where $\mathcal{Z}(G)$ denotes the center of the group $G$. In 1906, Fite \cite{Fite} proved that for finite nilpotent groups, this condition is sufficient for the existence of $\delta_{irr}(G)$. Then, in 1911, Burnside \cite{Burnside} provided an example of a semidirect product of $(C_{3}\times C_{3})\rtimes C_{2}$, where $C_{2}$ acts on $C_{3} \times C_{3}$ without non trivial fixed points, whose center is trivial (cyclic) and has no faithful irreducible representation. This shows that the condition may not be sufficient for an arbitrary finite group. However, when $\delta_{irr}(G)$ exists, we have the inequality
$\delta(G) \leq \delta_{irr}(G).$ Hence, $\delta_{irr}(G)$ is an upper bound for $\delta(G)$, if $\delta_{irr}(G)$ exists. Burnside \cite{Burnside} proved that if $G$ does not contain two distinct minimal normal subgroups whose orders are powers of the same
prime, then $G$ has a faithful irreducible representation. Later, in 1948, Kochendorffer generalized Fite's work \cite{Fite} and proved that if all Sylow subgroups of $G$ have
cyclic center, then $G$ has a faithful irreducible representation. In 1956, Zmud \cite{Zmud} proved that a finite group $G$ has a faithful complex representation with $k$ irreducible constituents if and only
if the socle of $G$ can be generated by $k$ elements as a normal subgroup.
  Szechtman provided a brief history of determining $\delta_{irr}(G)$ in \S2 of \cite{Sc}. 
  Further, $\delta(G) \leq c(G) \leq \mu(G),$ where $c(G)$ denotes the minimal degree of a faithful representation of $G$ by complex quasi-permutation matrices and $\mu(G)$ denotes the least positive integer $n$ such that $G$ is embedded in the symmetric group $S_{n}$.  These two quantities $c(G)$ and $\mu(G)$ have been studied extensively for $p$-groups of odd order, VZ $p$-groups, Camina $p$-groups, $p$-groups with cyclic center, and many other classes of groups (see e.g. \cite{MR3813700, MR1439603,MR1438628,MR1708003,MR1811393,PU,MR4552902}).   
\par In this article, we begin with listing several examples in \S\ref{examples}.  Then, in Theorem \ref{t2}, we determine $\delta(G)$ for certain $p$-groups. \begin{theoremalph}(Theorem \ref{t2}) If $G$ is a finite $p$-group with $\mathcal{Z}(G)$ cyclic, then $\delta(G)$ equals $\delta_{irr}(G)$. Further, if $G$ is a $p$-group, then $\delta(G)$ has been precisely computed for certain classes of $p$-groups, in particular, it covers all the $p$-groups of order up to $p^{5}$.  \end{theoremalph}
    In the next theorem, we compute $\delta(G)$ based on the set of degrees of all the irreducible characters of $G$, i.e., $\operatorname{cd}(G)$. 
    \begin{theoremalph}(Theorem \ref{t9})  Let $G$ be a finite nonabelian $p$-group with cyclic center. Then 
    $\delta(G)=p^{a},$ if $\operatorname{cd}(G)=\{1,p^{a}\}$ or $\operatorname{cd}(G)=\{1,p, p^{a}\}$, where $a \geq 1$ is an integer. Further, $\delta(G)=p,$ if and only if, $\operatorname{cd}(G)=\{1,p\}.$
\end{theoremalph}

Now, we provide information regarding representation dimension of nonabelian nilpotent groups with cyclic centers.
\begin{theoremalph}(Corollary \ref{c1}) Let $G=P_{1}\times P_{2} \times \cdots \times P_{n}$ be a finite nonabelian nilpotent group with cyclic center, where $P_{i}'s$ are Sylow subgroups. Then $\delta(G)$ equals \linebreak  $\delta_{irr}(P_{1})+ \delta_{irr}(P_{2})+\cdots + \delta_{irr}(P_{n})$.\end{theoremalph} In the subsequent theorems, we compute representation dimension for nonabelian groups of odd order having exactly two nonlinear irreducible characters of each degree and for groups whose all nonlinear irreducible characters have distinct degrees.
\begin{theoremalph}(Theorem \ref{t3}) 
Let $G$ be a finite nonabelian group of odd order with exactly two nonlinear irreducible characters of each degree. Then $\delta(G)$ equals $\delta_{irr}(G)$ and has been precisely computed, in case $G$ is an extraspecial $3$-group or Frobenius group of odd order $\frac{(p^{n}-1)p^n}{2}$ for some odd prime $p$ and some integer $n$, with cyclic kernel $K$ of order $p^{n}$ and cyclic complement. Otherwise, an upper bound on $\delta(G)$ has been obtained when $K$ is abelian but non cyclic. 
\end{theoremalph}
\begin{theoremalph}(Theorem \ref{t4}) 
If $G$ is a finite nonabelian group whose all nonlinear irreducible characters have distinct degrees, then $\delta(G)=\delta_{irr}(G)$ and $\delta(G)$ has been precisely computed. \end{theoremalph}
In Appendix~\ref{gapcode}, we present GAP codes to compute $\delta(G)$ and $\delta_{irr}(G)$.  


\section{Examples}\label{examples}
 
\noindent We begin with some examples to understand the representation dimension, $\delta(G)$, of a finite group $G$.

\begin{example}\label{exc2c2} 
If $G$ is a cyclic group, then $\delta(G)=\delta_{irr}(G)=1$. It is evident that for a finite group $G$, $\delta(G)=1$, if and only if, $G$ is cyclic. If $G = C_{2} \times C_{2}$, then $\delta(G)$ equals $2$, whereas $\delta_{irr}(G)$ does not exists,
which shows that understanding the character table of a group is not enough to determine its representation dimension. In fact, it can be seen that the representation dimension of a finite abelian group is the minimal number of generators of the group.
 \end{example} 
\begin{example}
If $G$ is a dihedral group, quaternion group, or $C_n\times C_m$, where $m$ is not coprime to $n$, then $\delta(G)=2$. Determining all finite groups with representation dimension $2$ is equivalent to determining all finite subgroups of $\GL(2,\mathbb{C})$. It is pertinent to mention that Brauer, Feit and their successors classified finite linear groups of degree up to $11$  over $\mathbb{C}$ (for more details, see the survey article  \cite{TZ} by Tiep and Zalesskii).  
\end{example} 

\par If $G$ is a finite simple group then every non-trivial character is faithful. Hence $\delta(G) = \delta_{irr}(G)$ for these groups. However, for arbitrary groups, $\delta(G) \neq \delta_{irr}(G),$ is evident from the following example.
\begin{example}\label{exa4s3}
For the group  $G=A_{4} \times D_{10}$, we have $\delta(A_{4} \times D_{10})=5$ and  $\delta_{irr}(A_{4} \times D_{10})=6$.  
\end{example}

\begin{example} 
From the character table of the alternating group $A_{4}$, we obtain that $\delta(A_{4})=3$. We remark that the representation dimension of $2.A_4$, the double cover of $A_4$, is $2$. This is because it is a subgroup of ${\rm SU}(2)$. Along similar lines, $\delta(A_{5})=3$, as $A_5$ is a subgroup of ${\rm SO}(3)$. 

Based on the work of Bessenrodt-Tong-Viet-Zhang \cite[Lemma 3.1]{BTZ}, $\delta(A_{n}) = n-1$, when $n \geq 15$. Using the GAP code given in Appendix~\ref{gapcode}, we obtain that $\delta(A_{n}),$ for $6 \leq n \leq 14$, is $n-1$. Hence, the representation dimension of $A_{n}$, for $n \geq 6$, is $n-1$.
\end{example}  
\begin{example} 
 Due to Result 1 of \cite{R}, we have that $\delta_{irr}(S_n) = n-1$, when $n \geq 5$. Indeed, the standard character of $S_{n}$ is faithful and has degree $n-1$. Consequently, $\delta(S_{n})=n-1$, for $n \geq 5$.
\end{example} 
\begin{example} 
The minimal degree of a  faithful irreducible character of  ${\rm GL}(2,q)$, ${\rm GL}(3,q)$ and ${\rm GL}(4,q)$ is $q-1$, $q^{2}+q$ and $(q+1)(q^{2}+1)$ respectively (see tables II, V and IX  in \cite{Steinberg}). 
Since these corresponds to the minimal degree of a nonlinear irreducible character in the character table of their respective groups, these are also the representation dimension of the corresponding groups. 
\end{example}
\begin{example} We have that $\delta({\rm SL}(2,q))=q-1$, when $q >2$ is even, and $\delta({\rm SL}(2,q))$ is even integer that is nearest to $\frac{q}{2}$, when $q$ is odd (see \cite{MR354889}). The representation dimension for finite groups of Lie type is well studied, see~\cite{Lubeck, TZ} for further details.
\end{example} 
\begin{example} 
Let $G$ be a finite group with a non-normal Sylow subgroup $P$ such that the intersection of two distinct conjugates of $P$ is trivial. Then, from \cite[Theorem 3.2]{BLM} , it follows that $\delta(G) > \sqrt{|P|}-1$. 
\end{example}
\begin{example} If $H$ and $K$ are isoclinic groups with $|H|=|K|$, then $\delta(H)$ and $\delta(K)$ may not be equal as we can see in Theorem \ref{t2}. However, if $H$ and $K$ are finite $p$-groups of the same order of nilpotency class $2$ with derived subgroups being cyclic and centers being isomorphic, then $\delta(H)=\delta(K)$, due to Theorem 1.3 and Lemma 3.5 of \cite{BMS}.   
\end{example}
The essential dimension $\operatorname{ed}(G)$ of a finite group $G$ over $\mathbb{C}$, is also related to the representation dimension of $G$. The notion of the essential dimension of $G$  was introduced by Buhler and Reichstein \cite{BR}.  It has been proved in \cite[Proposition 4.15]{BF} that $\operatorname{ed}(G) \leq \delta(G)$. One is interested to see that when $\operatorname{ed}(G)$ coincides with $\delta(G)$.
\begin{example} In Theorem 4.1 of \cite{KM}, Karpenko and Merkurjev proved that if $G$ is a finite $p$-group, then $\operatorname{ed}(G)$ coincides with $\delta(G)$. Recently, in \cite{kaur2024essential}, the essential dimension of finite groups of order $\leq 63$ has been computed.
\end{example}

\section{Representation dimension of $p$-groups}\label{sec3}

 \noindent  In the next theorem, we prove that for finite $p$-groups whose center $\mathcal{Z}(G)$ is cyclic, $\delta(G)$ equals $\delta_{irr}(G)$ and also determine $\delta(G)$ for certain $p$-groups. Note that, for part (4) of the next theorem, we depend on the classification of nonabelian groups of order $p^5$, where $p$ is an odd prime, via their isoclinism classes $\Phi_2, \Phi_3, \cdots, \Phi_{10}$, as given in \cite{james1}. Also, for part (5) of the next theorem, we use the classification of nonabelian groups of order $2^5$ via their isoclinism classes $\tau_{2},~ \tau_{3},\cdots,~\tau_{8}$, as given in \cite{Hall}, known as Hall-Senior families.
 \begin{theorem}\label{t2} 
The representation dimension $\delta(G)$ of certain $p$-groups $G$ is as follows: 
\begin{itemize}
\item [(1)] $\delta_{irr}(G)$, if $G$ is a finite $p$-group with $\mathcal{Z}(G)$ cyclic;\item [(2)] $p$, if $G$ is a nonabelian $p$-group of order $p^{3}$; 
\item [(3)] $p$, if $G$ is a nonabelian $p$-group of order $p^{4}$ with $\mathcal{Z}(G)$ cyclic, and is $p+1$ otherwise;
\item [(4)] if $G$ is a nonabelian $p$-group of order $p^{5}$, where $p$ is an odd prime, then the following cases arise, using the classification of groups via their isoclinism classes given in \cite{james1};
\begin{itemize} \item [(i)] if $G \in \Phi_{2}$, then $\delta(G)$ equals $2$ in case $\mathcal{Z}(G)$ is cyclic, $p+1$ if $\mathcal{Z}(G) \cong   C_{p^{2}} \times C_{p}$ and is $p+2$ if $\mathcal{Z}(G) \cong C_{p} \times C_{p} \times C_{p}$; \item [(ii)] $p$ if $G \in \Phi_{9}$; \item [(iii)] $p+1$ or $2p$ if $G \in \Phi_{i}$, where $i \in \{3,4,6\}$;
\item [(iv)] $p^{2}$ if $G \in \Phi_{i}$, where $ i \in \{5,7,8,10\}$. 
 \end{itemize}
 \item [(5)] if $G$ is a nonabelian $p$-group of order $32$, then the following cases arise:
 \begin{itemize}\item [(i)] if $G \in \tau_{2}$, then $\delta(G)$ equals $2$ in case $\mathcal{Z}(G)$ is cyclic; $3$ if $\mathcal{Z}(G) \cong C_{4}\times C_{2}$; $4$ if $\mathcal{Z}(G) \cong C_{2}\times C_{2}\times C_{2}$;
 \item [(ii)] if $G \in \tau_{3}$, then $\delta(G)$ equals $2$ in case $\mathcal{Z}(G)$ is cyclic; $3$ or $4$ if $\mathcal{Z}(G) \cong C_{2} \times C_{2}$; \item [(iii)] $4$ if $G \in \tau_{i}$, where $i \in \{ 4,5,6,7,8\}.$\end{itemize}
\item [(6)] $\sqrt{[G:\mathcal{Z}(G)]}$, if $G$ is a $p$-group with $\mathcal{Z}(G)$ cyclic and $|G'|=p$;  
\item [(7)] $\sqrt{[G:\mathcal{Z}(G)]}$, if $G$ is isoclinic to a semi-extraspecial $p$-group with $\mathcal{Z}(G)$ cyclic (in particular, $p^{r}$, if $G$ is an extraspecial $p$-group of order $p^{2r+1}$); \item [(8)] $[G: A],$ if $G$ is a nonabelian normally monomial $p$-group with $\mathcal{Z}(G)$ cyclic, where $A$ is an abelian normal subgroup of maximal order (in particular, groups of order $p^{5}$ with its center cyclic); 
\item [(9)] $\sqrt{[G:\mathcal{Z}(G)]}+\delta(\mathcal{Z}(G))-1$, if $G$ is a $p$-group of nilpotency class 2 with its derived subgroup cyclic (in particular, $\sqrt{[G:\mathcal{Z}(G)]}$, if $G$ is a $p$-group of nilpotency class 2 with $\mathcal{Z}(G)$ cyclic).
\end{itemize}
\end{theorem}

\vspace{-\topsep} \begin{proof}\begin{itemize}\item [(1)] Let $G$ be a finite $p$-group with cyclic center.
Let $\rho$ be a faithful representation of degree $\delta(G)$. Due to \cite[Lemma 3.5]{BMS}, 
it follows that $\rho$ can be written as a direct sum of exactly $r$ irreducible representations, where $r$ is the minimal number of generators of $\mathcal{Z}(G)$. Since $\mathcal{Z}(G)$ is cyclic, it follows that $r=1$. Consequently, $\rho$ is an irreducible representation of $G$. Hence $\delta(G)$ equals $\delta_{irr}(G)$.
\item [(2)] Suppose $G$ is a nonabelian $p$-group of order $p^{3}$.
Clearly, no linear character of $G$ is faithful, otherwise, $G$ is cyclic.  Further, the order of $\mathcal{Z}(G)$ must be $p$, since if the order of $\mathcal{Z}(G)$ is $p^{2}$, then $G$ is abelian, which is not possible. Consequently, by \cite[Theorem 2.32]{Isaacs}, $G$ has a faithful irreducible character. Also, it follows from the above part that $\delta(G)$ equals $\delta_{irr}(G)$. So, we need to determine $\delta_{irr}(G)$. As the quotient group $G/\mathcal{Z}(G)$ is abelian, we have that $G' \leq \mathcal{Z}(G)$. Since the order of $\mathcal{Z}(G)$ is $p$, it follows that $G'$ equals $\mathcal{Z}(G)$. Consequently, due to \cite[Corollary 2.30]{Isaacs}, it follows that $\delta_{irr}(G)=p$ and hence $\delta(G) =p$. 
\item [(3)] Let $G$ be a nonabelian $p$-group of order $p^{4}$, where $p$ is an odd prime, and let $\mathcal{Z}(G)$ be cyclic. 
Then $[G:\mathcal{Z}(G)]  \leq p^{3}$. Consequently, the degree of an arbitrary nonlinear irreducible character must be $p$, due to \cite[Corollary 2.30]{Isaacs}. Since $\mathcal{Z}(G)$ is cyclic, $\delta_{irr}(G)$ exists and $\delta(G) = \delta_{irr}(G)$. Hence $\delta(G)=p$. 
Suppose that $\mathcal{Z}(G)$ is non cyclic. Then the order of $\mathcal{Z}(G)$ is $p^{2}$. From \cite[Theorem 2.32]{Isaacs}, it follows that the group $G$ does not have a faithful irreducible character. Further, any nonlinear irreducible character of $G$ must be of degree $p$, due to Corollary 2.30 of \cite{Isaacs}. Suppose $\chi$ is a faithful character of $G$. From \cite[Lemma 2.3]{MR}, it follows that $\chi$ is a sum of two irreducible characters of $G$. Consequently, the degree of $\chi$ is either $p+1$ or $2p$. However, in view of Lemma 15 of \cite{CKR}, if $p$ is an odd prime, then $\operatorname{deg}(\chi) \leq p+1$; therefore, $\operatorname{deg}(\chi) = p+1$. \par For an arbitrary group $G$ of order $16$ with non cyclic center, one can verify from the character table of $G$ that $\delta(G)=3$. Hence, the result follows.
\item [(4)] \begin{enumerate} \item [(i)] Let $G \in \Phi_{2}$. From Lemma 5.1 of \cite{PDG}, it follows that $|\mathcal{Z}(G)|=p^{3}$, $|G'|=p$ and the set of all the irreducible character degrees of $G$, i.e. $\operatorname{cd}(G)$, is $\{1,p\}$. Clearly $\delta(G)=p$, if $\mathcal{Z}(G)$ is cyclic. Suppose $\mathcal{Z}(G)$ is non cyclic. Now, $G$ is a $p$-group of nilpotency class $2$ and its derived subgroup is cyclic. Therefore, in view of Theorem 1.3 of \cite{BMS}, it follows that $\delta(G)=p+1$, if $\mathcal{Z}(G) \cong C_{p^{2}}\times C_{p}$ and $\delta(G)=p+2$,  if $\mathcal{Z}(G) \cong C_{p}\times C_{p}\times C_{p}$. \item [(ii)] Let $G \in \Phi_{9}$. Then, in view of Lemma 5.7 of \cite{PDG}, $\mathcal{Z}(G)$ is cyclic and ${\rm{cd}}(G)=\{1,p\}$. Consequently $\delta(G)=p$.\item [(iii)] Let $G \in \Phi_{3}$.  Lemma 5.3 of \cite{PDG} yields that $\mathcal{Z}(G)$ is a non cyclic subgroup of order $p^{2}$ and $\operatorname{cd}(G)=\{1,p\}.$ 
Consequently, from Lemma 3.5 of \cite{BMS}, we conclude that if $\chi$ is a faithful character of $G$, then $\chi=\chi_{1}+\chi_{2}$ is a direct sum of two irreducible characters.  
Hence $\delta(G)$ is $p+1$ or $2p$. 
\par Similar argument holds when $G \in \Phi_{i}$, where $i \in \{4,6\},$ in view of Lemma 5.2 and Lemma 5.4 of \cite{PDG}.
\item [(iv)] Let $G \in \Phi_{5}$. Then $G$ is an extraspecial $p$-group and hence $\delta(G)=p^{2}$ \cite[Theorem 21.2.18]{Karpilovsky}. If $G \in \Phi_{7}$, then $|\mathcal{Z}(G)|=p$ and hence cyclic. Therefore, $\delta_{irr}(G)$ exists. Further, the  $p$-groups of order at most $p^{5}$ are metabelian and hence normally monomial (see \cite{How}). For the convenience of the reader, we recall that a finite group is normally monomial if all its irreducible complex characters are induced from linear characters of normal subgroups (see \cite{MR2228628}). It follows from \cite[Proposition 3]{Mann} that $\delta_{irr}(G)$ equals the maximum degree of all the irreducible characters of $G$ and hence equals $p^{2}$, due to Lemma 5.5 of \cite{PDG}. Consequently, (1) of this theorem yields that $\delta(G)=p^{2}.$ Similar argument holds when $G \in \Phi_{8}$. If $G \in \Phi_{10}$, then $\delta_{irr}(G)=p^{2}$, in view of Lemma 5.8 of \cite{PDG}. Hence $\delta(G)=p^{2}$. \end{enumerate}
\item [(5)] \begin{itemize} \item [(i)] Since $G$ is a $2$-group of nilpotency class $2$ with its derived subgroup cyclic, the result follows from Theorem 1.3 of \cite{BMS}.\item [(ii)] If $G \in \tau_{3}$ and its center is cyclic, then $\delta_{irr}(G)$ exists. However, $\operatorname{cd}(G)=\{1,2\}$ and therefore, $\delta(G)=\delta_{irr}(G)=2$.
\par Let $G \in \tau_{3}$ with $\mathcal{Z}(G) \cong C_{2} \times C_{2}$. In view of Theorem 2.2 of \cite{bioch1976n}, isoclinic groups have the same set of irreducible character degrees. Using GAP \cite{GAP}, it can be verified that ${\rm{cd}}(G)=\{1,2\}$.  
From Lemma 3.5 of \cite{BMS}, it follows that if $\chi$ is a faithful character of $G$, then $\chi=\chi_{1} + \chi_{2}$. 
Hence $\delta(G)$ is at most $4$. \item [(iii)] If $G \in \tau_{4}$ or $\tau_{5}$, then the nilpotency class of $G$ is $2$ and hence Theorem 1.2 of \cite{BMS} yields that $\delta(G)=4$. If $G \in \tau_{6}$, then either $G$ is SmallGroup$(32,6)$ or SmallGroup$(32,7)$. In either case, $\mathcal{Z}(G)$ is cyclic and it follows that $\delta(G)$ equals $2$. Similarly, one can see that if $G \in \tau_{7}$ or $\tau_{8}$, then $\delta(G)$ equals $2$.  \end{itemize}
\item [(6)]  Lemma 1 of \cite{MT} yields that all non linear irreducible characters of $G$ are faithful of degree $\sqrt{[G:\mathcal{Z}(G)]}$. Hence, the result follows.  
 \item [(7)] Before proceeding, recall that a $p$-group $G$ is semi-extraspecial if $G$ satisfies the property that for every maximal
 subgroup $N$ of $\mathcal{Z}(G)$, $G/N$ is an extraspecial group (see \cite{MR3922647}). Now, suppose $G$ is isoclinic to a semi-extraspecial $p$-group with $\mathcal{Z}(G)$ cyclic. It follows from Theorem A of \cite{FM} that the degree of an arbitrary nonlinear irreducible character is $\sqrt{[G:\mathcal{Z}(G)]}$. Since $\mathcal{Z}(G)$ is cyclic, $G$ has a faithful irreducible character. Therefore, $\delta(G) =\sqrt{[G:\mathcal{Z}(G)]}$.
\item [(8)] The group $G$ has a faithful irreducible character $\chi$, since $\mathcal{Z}(G)$ is cyclic. By \cite[Lemma 4]{MR2228628}, $\delta_{irr}(G)=[G:A]$, where $A$ is an abelian normal subgroup of $G$ of maximal order. Consequently, $\delta(G)=[G:A]$. 
\item [(9)] This is a direct consequence of  Theorem 1.3 of \cite{BMS}. \end{itemize}
 \end{proof} 
Let $G$ be a finite $p$-group. In the next theorem, we provide a connection between $\delta(G)$ and $\operatorname{cd}(G)$.

\begin{theorem} \label{t9} Let $G$ be a finite nonabelian $p$-group with cyclic center. 
\begin{itemize} 
\item [(i)] $\delta(G)=p$ if and only if, $\operatorname{cd}(G) = \{1,p\}$;
\item [(ii)] if for some integer $a > 1$, $\operatorname{cd}(G)=\{1,p^{a}\}$ or $\operatorname{cd}(G)=\{1,p, p^{a}\}$, then $\delta(G)=p^{a}$.
\end{itemize}\end{theorem}
\vspace{-\topsep} \begin{proof}\begin{itemize}\item [(i)] Let $\delta(G)=p$. In view of Theorem \ref{t2}(1), $\delta(G)$ equals $\delta_{irr}(G)$. Hence $\delta_{irr}(G)=p$. Theorem 22.25 of \cite{Berkovich} states that if $G$ is a $p$-group with a faithful irreducible character of degree $p^{n}$, then it has derived length at most $n+1$. Consequently, in this case, the derived length of $G$ is at most $2$. This implies that $G$ is a metabelian group. Since metabelian groups are normally monomial, it follows that the set of all possible irreducible character degrees is $\{1,p\},$ due to \cite[Proposition 3]{Mann}. \par Now, assume that $\operatorname{cd}(G)=\{1,p\}$. Since $G$ is a finite  nonabelian $p$-group with $\mathcal{Z}(G)$ cyclic, it follows from Theorem \ref{t2}(1) that $\delta_{irr}(G)$ exists and $\delta(G)=\delta_{irr}(G)$.  Consequently $\delta(G)=\delta_{irr}(G)=p$. 
\item [(ii)] If $\operatorname{cd}(G) = \{1,p^{a}\}$, then it immediately follows from Theorem \ref{t2}(1) that $\delta(G)=\delta_{irr}(G)=p^{a}$. \par Suppose
$\operatorname{cd}(G) = \{1, p, p^a\}$. In this case, Lemma 21 of \cite{PU} yields that $\delta(G)=p^{a}$.\end{itemize}
\end{proof}  
\section{Representation dimension of the direct product of groups}

\noindent In \S3, we have proved that if $G$ is a finite $p$-group with cyclic center, then $\delta(G) = \delta_{irr}(G)$. In 1906, Fite \cite{Fite} proved that for a finite nilpotent group, a necessary and sufficient condition for the existence of a faithful irreducible character is that its center is cyclic. 
In this section, we provide information on the representation dimension of some classes of direct product of groups, in particular, nilpotent groups.
\begin{theorem}\label{edproduct}
Let $G_{1}$ and $G_{2}$ be two finite groups and let $G= G_{1}\times G_{2}$. Then the representation dimension $\delta(G)$ equals $\delta_{irr}(G_{1})+\delta_{irr}(G_{2})$: 
\begin{itemize}
\item [(1)]  if $G_1$ and $G_2$ are non cyclic simple groups; 
\item [(2)]  if $G_1$ and $G_2$ are non cyclic groups of prime power order (distinct primes) and with cyclic center;
\item [(3)] if $G_1$ and $G_2$ are non abelian monolithic groups, i.e., groups with a unique minimal normal subgroup that is contained in every non trivial normal subgroup.
\end{itemize} 
\end{theorem} 
\begin{proof}   
(1) Let $G$ be a non cyclic simple group. If $G$ has a faithful linear character, then $G$ is cyclic of prime order. Hence $G$ has no faithful linear character and every irreducible nonlinear character is faithful. Therefore, $\delta(G)=\delta_{irr}(G)$. \par
Let $\rho$ be a faithful representation of degree $\delta(G)$ and $\chi$ be the faithful character corresponding to it. One can readily verify that $\chi$ must be written as $\chi_{1}+\chi_{2}$  with $\chi_{1}(1)=\delta_{irr}(G_{1})$ and $\chi_{2}(1)=\delta_{irr}(G_{2})$. Hence $\delta(G_{1} \times G_{2})$ equals $\delta_{irr}(G_{1})$ + $\delta_{irr}(G_{2}).$  \vspace{.2cm}\\  
(2) Due to Theorem \ref{t2}(1), $\delta(G_{i})=\delta_{irr}(G_{i})$, for $i =1,2$. If $\chi$ is a minimal faithful character of $G_{1} \times G_{2}$, then $\chi =\chi_{1} +\chi_{2}+\cdots +\chi_{n}$, where $\chi_{i} \in \operatorname{Irr}(G_{1} \times G_{2})$, for $1 \leq i \leq n$. Let $\chi_{i}=\tau_{i}\psi_{i}$, for $1 \leq i \leq n$, where $\tau_{i} \in \operatorname{Irr}(G_{1})$ and $\psi_{i} \in \operatorname{Irr}(G_{2})$.  Since $G_{1}$ and $G_{2}$ are of coprime order, $ \langle (g_{1},g_{2})~|~ g_{1} \in  \underset{1 \leq i \leq n} {\cap}\operatorname{ker}(\tau_{i}),
g_{2} \in \underset{1 \leq i \leq n}{\cap} \operatorname{ker}(\psi_{i}) \rangle \leq \underset{1 \leq i \leq n}{\cap} \operatorname{ker}(\tau_{i}\psi_{i})=\operatorname{ker}(\chi) =\{1\},$ where $\operatorname{ker}(\chi)$ denotes the kernel of the character $\chi$. This implies that $\langle (g_{1},g_{2})~|~ g_{1} \in  \underset{1 \leq i \leq n} {\cap}\operatorname{ker}(\tau_{i}),
g_{2} \in \underset{1 \leq i \leq n}{\cap} \operatorname{ker}(\psi_{i}) \rangle =\{1\}$. 
Since $\mathcal{Z}(G_{1})$ is cyclic and every minimal normal subgroup of finite $p$-group is elementary abelian and central, it follows that  minimal normal subgroup of $G_{1}$ is unique of order $p$ contained in $\mathcal{Z}(G_{1})$.  
Further, if $\operatorname{ker}(\tau_{i})$ is non trivial, then $\operatorname{ker}(\tau_{i}) \cap \mathcal{Z}(G_{1})$ is non trivial, and hence $\operatorname{ker}(\tau_{i})$ contains the unique minimal normal subgroup of $G_{1}$ of order $p$. Therefore, $\operatorname{ker}(\tau_{i})=\{1\}$, for some $i$, otherwise  $\underset{1 \leq i \leq n} {\cap}\operatorname{ker}(\tau_{i})$ would contain the unique minimal normal subgroup of $G_{1}$, contradicting that $\underset{1 \leq i \leq n} {\cap}\operatorname{ker}(\tau_{i})=\{1\}$.  Similar argument holds for $\operatorname{ker}(\psi_{j})$, for some $j$. Hence $\operatorname{ker}(\psi_{j})=\{1\}$, for some $j$. Since $\chi$ is the minimal faithful character of $G$, it follows that $\chi$ must be written as $\chi_{1}+\chi_{2}$  with $\chi_{1}(1)=\delta_{irr}(G_{1})$ and $\chi_{2}(1)=\delta_{irr}(G_{2})$. Hence $\delta(G_{1} \times G_{2})$ equals $\delta_{irr}(G_{1})$ + $\delta_{irr}(G_{2}).$  \vspace{.2cm}\\
(3) Let $G$ be a monolithic group and let $\chi$ be a faithful character of $G$ with $\deg(\chi)=\delta(G)$. We write $\chi=\underset{1 \leq i \leq n}{\sum}\psi_{i}$, where $\psi_{i}$ are irreducible constituents of $\chi$.  Now, $\{1\} = \operatorname{ker}(\chi)=\underset{1 \leq i \leq n}{\cap} \operatorname{ker}(\psi_{i})$. 
We observe that $\operatorname{ker}(\psi_{j}) = \{1\}$ for some $j$, else
$\operatorname{ker}(\chi)$ would contain the unique minimal normal subgroup of $G$, contradicting that $\chi$ is faithful. 
Thus, $\delta(G) \leq \deg(\psi_j) \leq \deg(\chi) = \delta(G).$
Hence $\chi$ is a faithful irreducible character of $G$. Consequently, in this case, $\delta(G_{i})=\delta_{irr}(G_{i})$, for $i=1,2$. \par If $\chi$ is a minimal faithful character of $G_{1} \times G_{2}$, then $\chi =\chi_{1} +\chi_{2}+\cdots +\chi_{n}$, where $\chi_{i} \in \operatorname{Irr}(G_{1} \times G_{2})$, for $1 \leq i \leq n$. Let $\chi_{i}=\tau_{i}\psi_{i}$, for $1 \leq i \leq n$, where $\tau_{i} \in \operatorname{Irr}(G_{1})$ and $\psi_{i} \in \operatorname{Irr}(G_{2})$.  Since $G_{1}$ and $G_{2}$ are of coprime order, $ \langle (g_{1},g_{2})~|~ g_{1} \in  \underset{1 \leq i \leq n} {\cap}\operatorname{ker}(\tau_{i}),
g_{2} \in \underset{1 \leq i \leq n}{\cap} \operatorname{ker}(\psi_{i}) \rangle \leq \underset{1 \leq i \leq n}{\cap} \operatorname{ker}(\tau_{i}\psi_{i})=\operatorname{ker}(\chi) =\{1\}$. This implies that $\langle (g_{1},g_{2})~|~ g_{1} \in  \underset{1 \leq i \leq n} {\cap}\operatorname{ker}(\tau_{i}),
g_{2} \in \underset{1 \leq i \leq n}{\cap} \operatorname{ker}(\psi_{i}) \rangle =\{1\}$. Since $\operatorname{ker}(\tau_{i})$ is a normal subgroup of $G_{1}$, it contains the unique minimal normal subgroup. 
Therefore, $\operatorname{ker}(\tau_{i})=\{1\}$, for some $i$, otherwise $\underset{1 \leq i \leq n} {\cap}\operatorname{ker}(\tau_{i})$ would contain the unique minimal normal subgroup of $G$, contradicting that $\underset{1 \leq i \leq n} {\cap}\operatorname{ker}(\tau_{i})=\{1\}$. Similar argument holds for $\operatorname{ker}(\psi_{j})$. Therefore, $\chi$ must be written as $\chi_{1}+\chi_{2}$  with $\chi_{1}(1)=\delta_{irr}(G_{1})$ and $\chi_{2}(1)=\delta_{irr}(G_{2})$. Hence $\delta(G_{1} \times G_{2})$ equals $\delta_{irr}(G_{1})$ + $\delta_{irr}(G_{2}).$ 
\end{proof} 
\noindent Due to Theorem \ref{edproduct}(2), the following corollary arises:  
\begin{corollary} \label{c1} Let $G=P_{1}\times P_{2} \times \cdots \times P_{n}$ be a finite nonabelian nilpotent group with cyclic center, where $P_{i}'s$ denote Sylow subgroups. Then $\delta(G)=\delta_{irr}(P_{1})+ \delta_{irr}(P_{2})+\cdots + \delta_{irr}(P_{n})$.\end{corollary}
\section{representation dimension of groups with few irreducible characters}\label{sec4}

Chillag and Herzog~\cite{CH} provided a complete classification of nonabelian groups of odd order that have exactly two nonlinear irreducible characters of each degree. We use this classification to determine the representation dimension of such groups.
\begin{theorem}\label{t3} 
Let $G$ be a finite nonabelian group of odd order with exactly two nonlinear irreducible characters of each degree. Then
\begin{itemize} 
\item[(i)]  $\delta(G) = 3^{r}$, if $G$ is an extraspecial $3$-group of order $3^{2r+1}$, for some integer $r \geq 1$. 
\item[(ii)]  $\delta(G) =\frac{p^{n}-1}{2}$, if $G$ is a Frobenius group of odd order $\frac{(p^{n}-1)p^n}{2}$ for some odd prime $p$ and some integer $n$, with cyclic kernel $K$ of order $p^{n}$ and cyclic complement. In case, $K$ is abelian but non cyclic, then $\delta(G)$ is atmost $\frac{(p^{n}-1)}{2}\operatorname{rank}(K)$.
\end{itemize} 
\end{theorem} 
\begin{proof}  
From Theorem 1 of \cite{CH}, a nonabelian group $G$ of odd order which has exactly two non linear irreducible characters of each degree, is one of the following: 
\begin{itemize} 
\item[(i)] an extraspecial $3$-group of order $3^{2r+1}$, for some integer $r \geq 1$; or
\item[(ii)] a Frobenius group of odd order $\frac{(p^{n}-1)p^n}{2}$ for some odd prime $p$ and some integer $n$, with abelian Frobenius kernel $K$ of order $p^{n}$ and cyclic Frobenius complement.
\end{itemize} 
Clearly, (i) follows from Theorem \ref{t2}(7). For (ii), first consider the case when the Frobenius kernel $K$ of the Frobenius group $G$ is cyclic, and so is its complement. We denote the Frobenius complement by $H$. It follows from Theorem 37.5.1 of~\cite{Karpilovsky} that the irreducible characters of $G$ are $\chi_{1}^{G},~\chi_{2}^{G},\cdots,~\chi_{r}^{G}, \lambda_{1},~\lambda_{2},\cdots,~\lambda_{s}$, where $\{\chi_{1},\chi_{2},\cdots,\chi_{r}\}$ is a complete set of representatives of $G$-conjugacy classes of non principal irreducible characters of $K$, $\chi_i^G$ are their induced representations on $G$, for $1 \leq i \leq r$, and $\lambda_{1},\lambda_{2},\cdots,\lambda_{s}$ are extensions of irreducible characters of $H$ to $G$. Clearly, $\lambda_{j}$'s are linear non faithful characters, for $1 \leq j \leq s$. Note that the degree of each nonlinear irreducible character of $G$ is $[G: K]$. Since $K$ is cyclic, $G$ has a nonlinear irreducible character which is faithful of degree $[G:K]$, i.e., $\frac{p^{n}-1}{2}$. As it is the minimal nonlinear irreducible character degree, it immediately follows that $\delta(G)= \delta_{irr}(G)=\frac{p^{n}-1}{2}$. 
\par Now suppose that the Frobenius kernel $K$ of $G$ is abelian but non cyclic and its complement $H$ is cyclic.  Note that $K \leq \operatorname{ker}(\lambda_{i})$, for $1 \leq i \leq s$. Therefore, $\delta(G)$ is atmost $[G:K]\delta(K)=\frac{(p^{n}-1)}{2}\operatorname{rank}(K)$. This completes the proof of the theorem. \end{proof} 

Theorem 2.3 of \cite{DNT} provides the classification of finite metabelian groups of odd order with the property that any two non principal irreducible characters of the same degree are Galois conjugate. Once again, we can compute the representation dimension for such groups using this classification. 
\begin{theorem} 
Let $G$ be a finite metabelian group of odd order with the property that any two non principal irreducible characters of the same degree are Galois conjugate. Then
\begin{itemize}
\item[(i)] $\delta(G)=\delta_{irr}(G)=p$, where  $[G:K] =p$ is a prime, if the group $G$ is a Frobenius group with Frobenius complement of order $p$ and Frobenius kernel $K$ is cyclic group of prime order $q$, where $q \neq p$.
\item[(ii)] $\delta(G)$ is atmost $\operatorname{rank}(K).p$, if  $G$ is a Frobenius group with complement of prime order $p$ and the Frobenius kernel $K$ is an elementary abelian $q$-group of order $q^{b}$, where $q \neq p$ is a prime and $b$ is also a prime with $p=\frac{q^{b}-1}{q-1}$.
\end{itemize}
\end{theorem} 

\section{representation dimension for groups in which the degrees of nonlinear irreducible characters are distinct}

Let $G$ be a finite group whose all nonlinear irreducible characters have distinct degrees. In this section, we compute representation dimension for such groups. These groups are classified by Berkovich, Chillag, and Herzog in~(\cite{BCH}, pg.955) and the explicit description is as follows.
\begin{theorem*}[Berkovich-Chillag-Herzog] 
If $G$ is a nonabelian finite group whose all \linebreak nonlinear irreducible characters have distinct degrees. Then, one of the following holds:
\begin{itemize}
\item[(1)] $G$ is an extraspecial $2$-group. 
\item[(2)] $G$ is a Frobenius group of order $p^{n}(p^{n} - 1)$ for some prime power $p^{n}$ with an abelian Frobenius kernel of order $p^{n}$ and a cyclic Frobenius complement. 
\item [(3)] $G$ is a Frobenius group of order $72$ in which the Frobenius complement is isomorphic to the quaternion group of order $8$. 
\end{itemize}
\end{theorem*}
\noindent Now, we have the following:
\begin{theorem}\label{t4} 
Let $G$ be a nonabelian finite group whose all nonlinear irreducible characters have distinct degrees. Then, $\delta(G)=\delta_{irr}(G)$, and it equals,
\begin{itemize} 
\item [(i)] $2^{r}$, if $G$ is an extraspecial $2$-group of order $2^{2r+1}$, for some integer $r \geq 1$. 
\item [(ii)] $p^{n}-1$, if $G$ is a Frobenius group of order $p^{n}(p^{n}-1)$ for some prime power $p^{n}$ with an abelian Frobenius kernel of order $p^{n}$ and a cyclic Frobenius complement. 
\item [(iii)] $8$, if $G$ is a Frobenius group of order $72$ in which the Frobenius complement is isomorphic to the quaternion group of order $8$. 
\end{itemize}
\end{theorem}
\begin{proof} Let $G$  be a nonabelian finite group whose all nonlinear irreducible characters have distinct degrees. Then, due to the above theorem, we can deal with such groups case by case. \vspace{.25cm}\\
\noindent{\bfseries Case 1.} $G$ is an extraspecial $2$-group. \\ In this case, $\delta(G)=\delta_{irr}(G)=2^{r}$, in view of Theorem~\ref{t2}(7).\vspace{.25cm}\\
\noindent{\bfseries Case 2.} $G$ is a Frobenius group of order $p^{n}(p^{n}-1)$ for some prime power $p^{n}$ with an abelian Frobenius kernel $K$ of order $p^{n}$ and a cyclic Frobenius complement $H$. \\ It follows from Theorem 5.1 of \cite{Karpilovsky} that all the irreducible characters of $G$ are obtained from the extensions of the irreducible characters of $H$ to $G$, and the rest irreducible characters of $G$ are induced from the representatives of $G$-conjugacy classes of non principal irreducible characters of $K$.
 As $K$ is abelian, the degree of each nonlinear irreducible character is $[G: K]$. By assumption, all nonlinear irreducible characters of the group $G$ have distinct degrees. This yields that there is a unique nonlinear irreducible character of degree $[G:K]$, which is faithful, due to Lemma 1.4 of \cite{IS}. Consequently, $\delta(G)=\delta_{irr}(G)=p^{n}-1$. \vspace{.25cm}\\
\noindent{\bfseries Case 3.} $G$ is a Frobenius group of order $72$ in which the Frobenius complement is isomorphic to the quaternion group of order $8$. \\ From Lemma 1.4 of \cite{IS}, it follows that $G$ has a unique nonlinear irreducible faithful character of degree 8. This completes the proof of the theorem. 
\end{proof} 

Iranmanesh and Saiedi \cite[Lemma 1.4]{IS} proved that if $G$ is a finite group with a unique nonlinear irreducible character, then it is faithful. In view of the above theorem, we have the following:
\begin{corollary} 
If $G$ is a finite group with a unique nonlinear irreducible character, then $\delta(G)=\delta_{irr}(G)$ is 
\item[(i)] $2^{r}$, if $G$ is an extraspecial $2$-group of order $2^{2r+1}$. 
\item[(ii)] $p^{n}-1$, if $G$ is a Frobenius group of order $p^{n}(p^{n}-1)$ for some prime power $p^{n}$ with an abelian Frobenius kernel of order $p^{n}$ and cyclic Frobenius complement. 
\end{corollary} 
\appendix

\section{GAP codes to compute $\delta(G)$ and $\delta_{irr}(G)$}\label{gapcode}
In this section, we present computational algorithms using GAP \cite{GAP} to determine the representation dimension and minimal faithful irreducible character degree of a finite group. \par Given a finite group $G$, we describe GAP code in Algorithm~\ref{A2} to compute the minimal faithful irreducible character degree, i.e., $\delta_{irr}(G)$. For this, we first give Algorithm~\ref{A1} to compute a list of all the irreducible characters of $G$ and their kernels. In Algorithm~\ref{A1}, we have a list $M$ of all the irreducible complex characters of $G$ with their values on conjugacy classes. Suppose there are $s$ irreducible characters in the list, i.e., $s$ number of conjugacy classes. For the $i$th irreducible character $M[i]$ mentioned in the list, look for that $ j \in [1 \cdots s]$ for which the value of $M[i]$ on $jth$ conjugacy class coincides with its degree, i.e., $M[i][j]=M[i][1]$. The union of all these conjugacy classes is the kernel $K$ of $M[i]$. Repeat this process for each irreducible character. Suppose $L=[~]$. In the list $L$, add $[i,K]$ for each $i \in [1 \cdots s]$. The list $L$ is the desired one. \par The first step of the Algorithm~\ref{A2} calls in the character table of $G$ and the number of conjugacy classes of $G$. Denote by $M$, the list of irreducible characters of $G$ with their values on conjugacy classes. 
From $L$, i.e., list of all irreducible characters of $G$ and their kernels, list those irreducible characters whose kernel is $\{e\}$ and denote the list by $F$. Consider the list $D$ of degrees of irreducible characters in $F$, i.e., the degrees of irreducible faithful characters of $G$. If $D$ is empty, then output is $0$, i.e., $\delta_{irr}(G)$ does not exist. If $D$ is nonempty, then the minimal irreducible faithful character degree $\delta_{irr}(G)$ is returned. 
 \par  
\begin{algorithm}
\begin{scriptsize}
\KwData{The list $M$ of all irreducible characters of finite group $G$ with their values on conjugacy classes.}
ListKernels:= function(M)\;
s := Size(M)\;
L := [~]\;
\For{i in \rm{[1..s]}}{K := [~]\;
{\For{j in \rm{[1..s]}}{\If{\rm{M[i][j]~=~M[i][1]}}{Add(K,j);}}Add(L,[i,K]);}}return L\\
\KwResult{Returns a list consisting of [$\chi, \operatorname{ker}(\chi)$], where $\chi$ is an irreducible character of $G$}
\caption{List of all the irreducible characters of $G$ and their kernels}\label{A1}
\end{scriptsize}
\end{algorithm}

\begin{algorithm}
\begin{scriptsize}
\KwData{A finite group $G$.}
c:=CharacterTable(G)\;
CC := ConjugacyClasses(G)\;
s := Size(CC)\;
M := List(Irr(c),ValuesOfClassFunction)\;
L := ListKernels(M)\;
{F := \rm{Filtered}(L, x $->$ \rm{Length}(x[2]) = 1)\;}{D := [~]\;}{d := 0\;}\For{j in \rm{F}}{Add(D,M[j[1]][1]);}\If{\rm{not(IsEmpty(D))}}{d := Minimum(D);} return d\\
\KwResult{$\delta_{irr}(G)$ }
\caption{Minimal faithful irreducible character degree of $G$}\label{A2}
\end{scriptsize}
\end{algorithm}

\noindent Next, we provide Algorithm~\ref{A8} to determine the representation dimension of $G$. If $G$ is a finite abelian group, then the representation dimension of $G$ is the minimal number of generators of the group $G$, in view of Example 2.1. Hence, the first step is to provide Algorithm~\ref{A5} for determining the representation dimension of a finite abelian group $G$. Let $G$ be a finite abelian group. We determine the structure of $G$ as a direct product of cyclic groups of prime power order, using the GAP command AbelianInvariants($G$). Next, we provide an algorithm to convert AbelianInvariants list into minimal cyclic factors. In this, we need an algorithm to delete all the occurences of $r$ in a list $S$. In the Algorithm~\ref{A8}, we also add the information that has been obtained from Theorem \ref{t2}. After this, the next step is to provide Algorithm~\ref{A6} for determining all the possible kernels of characters of $G$ and Algorithm~\ref{A7} to determine the degree of each character of $G$.

\begin{algorithm}
\begin{scriptsize}
  \KwData{A finite list $S$.}
  DelAll := function(S,r)\;
  \While{\rm{Position(S,r) $<>$ fail}}{Remove(S, Position(S,r))\;}{return S\;}
  \KwResult{Deletes all occurences of r in a list S.}
  \end{scriptsize}
  \caption{}\label{A3}
\end{algorithm}

\begin{algorithm}
\begin{scriptsize}
 \KwData{A finite list $T$.}
    CyclicBunches := function(T)\;
            S := [~]\;
        R := [~]\;
        \For{\rm{x in T} }{Add(S,x)}
                n := Size(S)\;
 \While{\rm{n $>$ 0}}{{{a := S[1]\;}{S[1] := 1\;}}
 \For{\rm{i in [2..n]}}{{\If{\rm{Gcd(a,S[i]) = 1}\;}{{a := a*S[i]\;}{S[i] := 1\;}}}}}
 {Add(R,a)\;}
 {DelAll(S,1)\;}
 {n := Size(S)\;}
        {return R\;}
\KwResult{If $T$ equals AbelianInvariants list, then result is the minimal cyclic factors.}
\end{scriptsize}
\caption{}\label{A4}
\end{algorithm}

\begin{algorithm}
    \begin{scriptsize}
        \KwData{A finite abelian group $G$}
        EmbeddingDegreeAbelian := function(G)\;
        T := AbelianInvariants(G)\;
        return Size(CyclicBunches(T))\;
    \KwResult{$\delta(G)$}
    \caption{Representation dimension of an abelian group $G$}\label{A5}
    \end{scriptsize}
\end{algorithm}

\begin{algorithm}
\begin{scriptsize}
\KwData{A finite group $G$.}
KernelIntersection := function(L,S)\;
$s$ := Size(L)\;
IntKer:= [1..s]\;
\For{i in S}
{IntKer := Intersection(IntKer, L[i][2]);}\
return IntKer;\\ 
\KwResult{Returns the intersection of $\operatorname{ker}(\chi)$ of characters $\chi$, where i lies in the set S.}
\end{scriptsize}
\caption{The set of all possible kernels of characters of a group $G$}\label{A6}
\end{algorithm}

\begin{algorithm}
\begin{scriptsize}
\KwData{A finite group $G$.}
DegreeSum := function(M,S)\;
d := 0\;
\For{\rm{i in S}}
{d:=d+M[i][1]\;}
return d;\\ 
\KwResult{Returns the sum of degrees of characters $\chi_{i}$, \rm{where~i~lies~in~the~set~S}.}
\end{scriptsize}
\caption{Degrees of all possible characters of group $G$}\label{A7}
\end{algorithm}
\newpage In the next algorithm, we use previous algorithms to determine the representation dimension of a group $G$.  
Clearly, the representation dimension of group $G$ is less than or equal to the degree of the regular character of $G$, denoted by MinDeg. Now, we look at the possible combination of irreducible characters that yield a faithful character whose degree is less than MinDeg. Suppose a given combination is providing a faithful character $\chi$ with its degree being less than MinDeg. Then, repeat the process with MinDeg being replaced by $\operatorname{deg}(\chi)$. This process will stop at a finite number of steps. Finally, the representation dimension of $G$ is determined.

\begin{algorithm}
\begin{tiny}
\KwData{A finite group $G$.}
RepresentationDimension:= function(G)\;
factors := Collected(Factors(Order(G)))\;
ZG := Center(G)\;
ZCyclic := IsCyclic(ZG)\;
GAbelian := IsAbelian(G)\;
\If{\rm{GAbelian}}{return EmbeddingDegreeAbelian(G)\;}
\If{\rm{Size(factors) = 1}}{{p := factors[1][1]\;}{expnt := factors[1][2]\;}\If{\rm{ZCyclic}}{return(IrredFaithfulMinDegree(G))\;}\If{\rm{expnt = 3 and not(GAbelian)}}{return factors[1][1]\;}\If{\rm{expnt = 4 and not(GAbelian)}}{\If{\rm{ZCyclic}}{{return p\;}\Else{return p+1\;}}}}\If{\rm{NilpotencyClassOfGroup(G) = 2~and~IsCyclic(DerivedSubgroup(G))}}{return Sqrt(Index(G,ZG))+EmbeddingDegreeAbelian(ZG)-1\;}
c := CharacterTable(G)\;
CC:= Conjugacy classes(G)\;
s:=Size(CC)\;
M := List(Irr(c),ValuesOfClassFunction)\;
L := ListKernels(M)\;
minDeg := Sum(TransposedMat(M)[1])\;
minK := [1..s]\;
rList := [1..s-1]\;
r := 1\;
\While{\rm{r $\leq$ minDeg}}{{rCombos := Combinations([1..s],~r)\;\For{\rm{x in rCombos}}{DegSumx := DegreeSum(M,x)\;{\If{\rm{DegSumx $<$ minDeg and r $<$ minDeg}}{\If{\rm{Length(KernelIntersection(L,x)) = 1}}{minDeg := DegSumx\;minK := x\;}}}}}{ r := r+1\;}}{return minDeg;\\}
\KwResult{$\delta(G)$}
\caption{Representation dimension of a finite group $G$}\label{A8}
\end{tiny}
\end{algorithm}

\newpage \subsection*{Acknowledgment} We thank B. Sury for his feedback on this paper which helped us improve it. The first named author acknowledges the research support of the Department of Science and Technology (INSPIRE Faculty No. DST/INSPIRE/04/2023/001200), Government of India. The third named author is funded by NBHM through \linebreak  02011/23/2023/NBHM(RP)/RDII/5955 for this research.


\bibliographystyle{amsplain}
		\bibliography{BIBTEX}
	\end{document}